\documentclass{mcom-l}

\newcommand{\R}{\mathbb{R}}
 
\newcommand{\Z}{\mathbb{Z}}
\newcommand{\N}{\mathbb{N}}

\newcommand{\Oe}{\tilde{O}}

\DeclareMathOperator{\ord}{ord}

\DeclareMathOperator{\polylog}{polylog}

\usepackage{amssymb}
\usepackage{graphicx}
\usepackage{xcolor}
\usepackage{booktabs}
\usepackage{tabularx}
\usepackage{adjustbox}
\usepackage{hyperref}
\usepackage{cleveref}
\usepackage{algpseudocode}
\usepackage{float}
\usepackage{caption}
\usepackage{subcaption}
\usepackage{enumerate}
\usepackage{enumitem}   
\usepackage{alltt}
\algtext*{EndWhile}
\algtext*{EndFor}
\algtext*{EndIf}

\copyrightinfo{}{}

\newtheorem{theorem}{Theorem}[section]
\newtheorem{lemma}[theorem]{Lemma}
\newtheorem{corollary}[theorem]{Corollary}

\theoremstyle{definition}
\newtheorem{definition}[theorem]{Definition}

\newtheorem{remark}[theorem]{Remark}
\newtheorem{assumption}[theorem]{Assumption}

\numberwithin{equation}{section}

\begin{document}

\title[Private key sieves for RSA]{Private key sieves for improved attacks on RSA}

%    author one information
\author[M. Hittmeir]{Markus Hittmeir}
\address{NORCE Research}
\curraddr{Nygårdsgaten 112, 5008 Bergen, Norway}
\email{mahi@norceresearch.no}
\thanks{The author was supported by the Research Council of Norway (grant 357539).}
\begingroup \makeatletter \renewcommand{\@makefntext}[1]{#1} 
\begin{NoHyper} 
	\footnotetext{
		\hangindent=1.5em 
		\hangafter=1 
		This is a short version of the paper. Implementation details will be included in a future version.
	} 
\end{NoHyper} 
\makeatother 
\endgroup

%    \subjclass is required.
\date{}
\dedicatory{}

%    Abstract is required.
\begin{abstract}
The RSA cryptosystem has been studied extensively and its insecurity in specific scenarios is well known. For instance, a sufficiently small private key allows to factorize the modulus $N$ and break the system. Attacks like this are often based on lattice reduction for finding small roots of polynomials, but high lattice dimensions can make it difficult to achieve the theoretically proved bounds in practice. In this paper, we adapt a sieve technique originally developed for deterministic integer factorization and apply it to RSA parameters. The sieve allows for an asymptotically significant reduction of the private key search space, and can be used to decrease the bounds for the polynomial roots. We also provide a square-root speed-up of the best currently known runtime bound for factoring semiprime $N=pq$ with small prime difference $q-p$. 
\end{abstract}

\maketitle

\section{Introduction}

RSA is a widely used cryptosystem that enables digital signatures, authentication and secure communication. It was first introduced by Rivest, Shamir and Adleman (\cite{RSA-rsa}) in 1978. The system is defined by a modulus $N=pq$, where $p$ and $q$ are large primes, a public exponent $e$ used for encryption, and a private key $d$ used for decryption. Let $\phi(N)=(p-1)(q-1)$ denote the value of Euler's totient function at $N$, then $e$ and $d$ satisfy $ed\equiv 1 \mod \phi(N)$. If one were able to compute the primes $p$ and $q$ from the public modulus $N$, one could compute $\phi(N)$, obtain $d$ as the inverse of $e$ modulo $\phi(N)$ and consequently break RSA. Therefore, the difficulty of finding the prime factors of natural numbers $N$, known as integer factorization, is a crucial assumption in the security of RSA.

The fastest currently known algorithms for integer factorization achieve runtime complexities that are subexponential in the bit-length of $N$. The two most effective general-purpose methods are the Quadratic Sieve, which works best for numbers up to around 300 to 350 bits, and the Number Field Sieve, which is the asymptotically fastest known practical factorization method. The monographs \cite{Rie-factorization} and \cite{Wag-joy} provide an overview on the variety of methods. However, the mentioned methods are not fast enough to break RSA systems using 2048-bit moduli, which is a commonly recommended key size (\cite{NIST}). Shor's algorithm (\cite{Sho-factoring}) can factorize integers in polynomial-time, but relies on quantum computing. Using a sufficiently large and fault-tolerant quantum computer, it would potentially be feasible to break RSA and similar cryptographic systems. While such devices do not exist yet, this prospect has motivated the development of post-quantum cryptography, aiming to construct systems not vulnerable to quantum attacks.

Certain mistakes in the choice of the prime factors $p$ and $q$ can make RSA vulnerable to non-quantum, special-purpose factorization methods. For instance, this is the case if $p$ and $q$ are very close together, or if either $p-1$ or $q-1$ is the product of only small primes. But if $p$ and $q$ are chosen carefully and satisfy certain properties, it is almost impossible to factorize 2048-bit moduli $N$ with any of the currently available techniques. 

There is a large variety of attacks that exploit other potential mistakes in the setup of the RSA cryptosystem. The reader may consult Boneh's survey (\cite{B-survey}) for an introductory and comprehensive overview of various methods and scenarios. In this paper, we are particularly interested in cases where either the public exponent $e$ or the private key $d$ is small, and where parts of $d$ have been exposed. These are classical attack scenarios that have been investigated in great detail. In 1998, Boneh, Durfee and Frankel (\cite{BDF-rsa}) showed that, if $e$ is very small such that a cost of $O(e)$ is feasible, then knowing a quarter of the least significant bits of $d$ allows to compute all of $d$. Blömer and May (\cite{BM-partialkeyexposure}) generalized this result to $e<N^{7/8}$, but their method requires a larger amount of known bits of $d$. The core idea of their result is to use the known bits of $d$ to construct a bivariate polynomial with a small root, using the already mentioned central equation $ed - 1 = k\cdot \phi(N)$, where $k\in\N$. The root is then of the shape $(k, p+q-1)$ and can be found in heuristic polynomial-time via lattice reduction techniques. Knowing the value of $p+q$ allows to factorize $N$ immediately.

In the case of small private keys $d$, Wiener's attack (\cite{W-smalld}) from 1990 shows that RSA is vulnerable as soon as $d<N^{1/4}$, since $d$ can then be recovered from the continued fraction expansion of $e/N$. Using lattice methods, Boneh and Durfee (\cite{BD-smalld}) improved this to the heuristic bound $d<N^{0.292}$. Ernst, Jochemsz, May and de Weger (\cite{EJMW-partialkeyexposure}) then developed a technique that allows to use knowledge of some bits of $d$ to relax the upper bound even further. Their approach can also be described as a problem of finding the roots of multivariate polynomials, and is thus similar to the method by Blömer and May described above.

\textbf{Our contributions.} We generalize a technique that has previously been studied in the context of deterministic integer factorization, called the hyperbolic sieve. It has first been developed in \cite{Hit-BSGS} and further investigated in \cite{Hit-subsetsum} and \cite{SMH-sieve}. The main idea is to consider sets of the shape 
\[
\mathcal{S}_{r}:=\{(x+x^{-1}N)\pmod r: x\in\Z_r^*\},
\]
where $r$ is a small prime. One can show that $\mathcal{S}_{r}$ contains only about $r/2$ of all possible residues modulo $r$, and one of those must be equal to $(p+q) \pmod r$. Considering a product $P$ of many small primes, we can use the Chinese Remainder Theorem to obtain a subset $\mathcal{S}_{P}$ of $\Z_P$ that is significantly smaller in size and contains the value $(p+q) \pmod P$. Such information can then be used in the search for $p+q$, e.g. as solution to the discrete logarithm problem
$
a^X\equiv a^{N+1}\mod N,
$
where $a\in\Z_N^*$ (see \cite{Hit-BSGS}). 

The generalization will allow us to apply the hyperbolic sieve not only to $p+q$, but also to other private RSA parameters such as $\phi(N)$ and $d$. We will modify the polynomials used in the attacks by Blömer and May (\cite{BM-partialkeyexposure}) and by Ernst, Jochemsz, May and de Weger (\cite{EJMW-partialkeyexposure}) to combine knowledge of $\phi(N)\pmod P$ with known bits of $d$. Our strategy is to loop through all candidates for $\phi(N)\pmod P$ obtained from the sieve and apply lattice reduction to try and find the roots of the resulting polynomials. We will succeed as soon as we use the correct value for $\phi(N)\pmod P$. This strategy is similar to the common approach of enumerating or guessing a few bits of unknowns like $p$, $p+q$ or $d$ to slightly reduce the bounds for the polynomial roots and the lattice dimension, thereby increasing the chance of success in practice (e.g., see \cite{WTS-enumeration}). However, the hyperbolic sieve allows for an \emph{asymptotically significant reduction} of the necessary number of guesses by a factor that is almost exponential in the size of the total search space. We thus achieve an extension of the attacks in \cite{BM-partialkeyexposure} and \cite{EJMW-partialkeyexposure} that is of both theoretical and practical interest.

In addition to the lattice-based attacks mentioned above, we also provide an improvement to an attack scenario that exploits small prime difference of the factors of $N$. It is well known that one can use Fermat's classical method to factor $N$ in around $\Lambda:=(q-p)^2/(4\sqrt{N})$ steps. Using the hyperbolic sieve, this has already been improved to around $\Lambda\cdot \exp(-C\log\Lambda/\log\log\Lambda)$ steps, where $C:=(1+o(1))\log 2$ (Theorem 4.6 in \cite{Hit-subsetsum}). We will now prove that $N$ can be factored in around
\[
\Lambda^{1/2}\cdot \exp\left(-\frac{C}{2}\log\Lambda/\log\log\Lambda\right)
\] 
steps. The algorithm is deterministic and the proof is fully rigorous.

\textbf{Overview.} \Cref{sec:prelim} briefly introduces lattices and well-known results concerning lattice reduction techniques. In \Cref{sec:sieve} we will show how to apply the hyperbolic sieve to RSA parameters. Section \ref{sec:smalle}, \ref{sec:smalld} and \ref{sec:smalldiff} concern our main results on RSA with small public exponent $e$, small private key $d$ and small prime difference, respectively.

\textbf{Notation.} Throughout this paper, $N$ always denotes an RSA modulus, i.e., a balanced semiprime $N=pq$ with two distinct prime factors $p$ and $q$, where $p<q<2\sqrt{N}$. 
By $\phi$ we denote the Euler's totient function. The pair $(N,e)$ denotes an RSA public key, and $d$ denotes an RSA private key. We write $\Z_n$ for the residue class ring $\Z/n\Z$ with $n\in\Z$, and $\Oe(\cdot)$ for an asymptotic bound suppressing polylogarithmic
factors. Concretely, $f\in O(\polylog g)$ means $f\in O(\log ^c g)$ for some constant $c\geq 0$, and $f\in\Oe(g)$ means $f\in O(g \polylog g)$.

\section{Lattice Reduction} \label{sec:prelim}
Let $\mathbf{b}_1,\ldots,\mathbf{b}_n \in \mathbb{R}^n$ be linearly independent vectors. The lattice generated by these vectors is defined as
\[
L = L(\mathbf{b}_1,\ldots,\mathbf{b}_n)
:=
\left\{
\sum_{i=1}^n a_i \mathbf{b}_i : a_i \in \mathbb{Z}
\right\}.
\]
The vectors $\mathbf{b}_1,\ldots,\mathbf{b}_n$ form a basis of $L$. If $B$ is the corresponding basis matrix, then the determinant of $L$ is defined as
$
\det(L):=|\det(B)|.
$
We call $n$ the dimension of the lattice. 

By $\|\cdot\|$ we denote the Euclidean norm.  For $\mathbf{b}=(b_1,\ldots,b_n)\in\Z^n$ and a polynomial $f(y,z)=\sum_{i,j}f_{i,j}y^iz^j\in\Z[y,z]$, we set 
\[ \|\mathbf{b}\|:=\Bigl(\sum_i b_i^2\Bigr)^{1/2}, \qquad \|f\|:=\Bigl(\sum_{i,j}|f_{i,j}|^2\Bigr)^{1/2}. \]

The following result concerns the well-known LLL lattice reduction algorithm by Lenstra, Lenstra and Lovász (\cite{LLL-lll}). The method's ability to find short vectors in lattices has a wide applicability.

\begin{lemma}[LLL]\label{thm:lll}
	Let $L\subseteq \Z^n$ be a lattice of dimension $n\geq 2$. In polynomial time, we may find two linearly independent vectors $\mathbf{v}_1,\mathbf{v}_2\in L$ that satisfy
	\[
	\|\mathbf{v}_1\|\leq \|\mathbf{v}_2\|
	\leq
	2^{n/4}\det(L)^{1/(n-1)}.
	\]
\end{lemma}

A standard strategy in lattice-based attacks on RSA is to construct multivariate polynomials $f$ with a small root modulo some natural number $n$. If certain conditions are satisfied,  Howgrave-Graham's (\cite{HG-smallroots}) variant of a result by Coppersmith \cite{C-smallsols} allows to deduce that $f$ has the same root as a polynomial over $\Z$. We will state this result for the case of bivariate polynomials. A similar formulation can be found in \cite[Lemma 1]{C-smallroots}.

\begin{lemma}[Howgrave-Graham] \label{thm:hg}
	Let $f(y,z) \in \Z[y,z]$ be a polynomial which is a sum of at most $\omega$ monomials. Suppose that
	$
	f(y_0,z_0) \equiv 0 \pmod{n}
	$
	such that $|y_0| \leq Y, |z_0| \leq Z$ and
	\[
	\|f(yY,zZ)\|
	<
	\frac{n}{\sqrt{\omega}}.
	\]
	Then $f(y_0,z_0)=0$ holds over the integers.
\end{lemma}

While we will use lattice reduction methods mostly as a black box in this paper, the later sections provide a more detailed sketch of how these results are applied in the framework used in \cite{BM-partialkeyexposure} and \cite{EJMW-partialkeyexposure}.

\section{The hyperbolic sieve for RSA parameters} \label{sec:sieve}
Let us now recall the definition of the hyperbolic sieve and the results that we will need in later sections. In addition, we will show how to extend the sieve to RSA parameters such as Euler's totient function $\phi(N)$ and the private key $d$.

{\definition{Let $N,P\in\N$ such that $\gcd(N,P)=1$. We define the sets
		\[
		\mathcal{H}_{P}:=\{(x,y)\in\Z_P^2 :xy\equiv N \mod P\}
		\]
		and
		\[
		\mathcal{S}_{P}:=\{x+y\mod P: (x,y)\in\mathcal{H}_{P}\}.
		\]}\label{def:hyper}
}

The set $\mathcal{H}_{P}$ is called a modular hyperbola and has the following properties:
\begin{enumerate}
	\item{$|\mathcal{H}_{P}|=\phi(P)$.}
	\item{If $N=uv$, then $(u\mod P, v \mod P)\in \mathcal{H}_{P}$.}
\end{enumerate}
Clearly, $N=uv$ implies that $(u+v)\pmod P\in\mathcal{S}_{P}$. One easily observes that $\mathcal{S}_{P}=\{(x+x^{-1}N)\pmod P: x\in\Z_P^*\}$. The set $\mathcal{S}_{P}$ has been introduced in \cite{Hit-BSGS}, where it has been used to search for the sum of the prime factors of semiprime numbers $N$. 

{\lemma{
		Let $N\in\N$ and $r$ be an odd prime with $r\nmid N$. Then
		\[
		|\mathcal{S}_{r}|=(r+(N|r))/2,
		\]
		where $(N|r)$ denotes the Legendre symbol.
}}

\begin{proof}
	This is Lemma 5.3 in \cite{Hit-BSGS}.
\end{proof}

{\theorem{
		Let $N\in\N$, $B\in\R^+$ and set $P=\prod_{2<r\leq B} r$ for primes $r$. Assume that $\gcd(N,P)=1$. Then
		\[
		|\mathcal{S}_{P}|\in \widetilde{O}\left(P\cdot \exp\left(-C\log P/\log\log P\right)\right),
		\]
		where $C=(1+o(1))\log 2$. The elements of $\mathcal{S}_P$ can be computed sequentially (requiring negligible space) and in asymptotically optimal time $O(|\mathcal{S}_P|\log P)$.
	}\label{thm:sievebound}}
\begin{proof}
	Let us first prove the statement about the cardinality of $\mathcal{S}_P$. We start by discussing the size of $B$. For Chebyshev's theta function $\vartheta(x):=\sum_{r\leq x} \log r$, it follows that
	\[
	e^{\vartheta(B)}=\prod_{r\leq B} r =2P.
	\]
	From the Prime Number Theorem, it follows that $\vartheta(x)-x\in o(x)$ (see \cite[p.80]{RR-primes}). We conclude that $B=(1+o(1))\log P$. According to Theorem 5.5 in \cite{Hit-BSGS}, we have
	\[
	|\mathcal{S}_{P}|\in O\left(\frac{P\log B}{2^{\pi(B)}}\right)\subseteq \widetilde{O}\left(\frac{P}{2^{\pi(B)}}\right).
	\]
	The Prime Number Theorem states that $\pi(B)\approx B/\log B$. Indeed, we have the explicit lower bound of $B/\log B<\pi(B)$ for $B\geq 17$ (see \cite[p.69]{RR-primes}). Hence, for $C=(1+o(1))\log 2$, we obtain
	\begin{align*}
		O(2^{-\pi(B)})= O\left(2^{-\frac{B}{\log B}}\right)=O\left(2^{-\frac{(1+o(1))\log P}{\log\left(\log P\right)}}\right)
		= O\left(\exp(-C\log P/\log\log P)\right),
	\end{align*}
	which proves the claim.
	
	The statement about computing the elements of $\mathcal{S}_P$ follows from Algorithm 4.1 and Theorem 4.2 in \cite{Hit-BSGS}, whose assumption $|\mathcal{S}_P|\geq \sum_{2<r\leq B} |\mathcal{S}_r|\cdot \log P$ is satisfied for sufficiently large $P$ (namely, for $B\geq 13$).
\end{proof}

Let $e$ be a public RSA exponent and $d$ be the private key. As mentioned in the introduction, we have 
\begin{equation}\label{eq:main}
	ed  - 1 = k\cdot \phi(N),
\end{equation}
where $k\in\N$.

{\label{lem:computek}\lemma{The value $k$ in (\ref{eq:main}) is an element of the set 
		\[
		\mathcal{K}:=\{(s-N-1)^{-1}\pmod e : s \in \mathcal{S}_e \wedge \gcd (e, N+1-s)=1\}.
		\]
	In particular, $|\mathcal{K}|\leq |\mathcal{S}_e|$.
	}}
\begin{proof}
	Note that $d=e^{-1}\pmod{\phi(N)}$. Hence, $k=(ed-1)/\phi(N)<ed/\phi(N)<e$, and we may write $k=k\pmod e$.
	In addition, we have $-1\equiv k \cdot \phi(N)\mod e$. Since $\gcd(e,\phi(N))=1$, we get
	\[
	k\equiv -\phi(N)^{-1}\equiv - (N+1-p-q)^{-1} \equiv (p+q-N-1)^{-1} \mod e,
	\]
	which implies that $k=k\pmod e\in\mathcal{K}$.
\end{proof}

{\lemma{\label{lem:phi} Let $P\in\N$ such that $\gcd(N,P)=1$. The value $\phi(N)\pmod P$ is an element of the set 
	\[
	\Phi_P:=\{(N+1-s)\pmod P : s \in \mathcal{S}_P\}.
	\]
	In particular, $|\Phi_P|\leq |\mathcal{S}_P|$.
}}
\begin{proof}
The claim follows immediately from the formula of Euler's totient function for semiprime numbers $N=pq$, namely $\phi(N)=(p-1)(q-1)=N+1-(p+q)$.
\end{proof}

{\label{lem:computed}\lemma{Let $P\in\N$ such that $\gcd(Ne,P)=1$. The value $d\pmod P$ is an element of the set
	\[
	\mathcal{D}_{P}(k):=\{e^{-1}(k(N+1-s)+1)\pmod P : s \in \mathcal{S}_P\}.
	\]
	In particular, $|\mathcal{D}_P(k)|\leq |\mathcal{S}_P|$.
}}
\begin{proof}
From (\ref{eq:main}) and $\gcd(e,P)=1$ it follows that 
\[
d\equiv e^{-1}(k\phi(N)+1)\equiv  e^{-1}(k(N+1-(p+q))+1)\mod P,
\]
which implies the claim.
\end{proof}

The sieve set for the private key $d$ depends on $k$, which is usually not known. However, if $e$ is small, then $|\bigcup_{k\in\mathcal{K}} \mathcal{D}_P(k)|$ may still be much smaller than $P$.

\section{Attack on RSA with small public exponent}\label{sec:smalle}
In 1998, Boneh, Durfee and Frankel \cite{BDF-rsa} showed that RSA is vulnerable to partial exposure of the private key. If $e$ is very small (such that enumeration over all candidates for $k<e$ is feasible), they proved that a quarter of the least significant bits of $d$ is enough to recover the complete key. In the same paper, they also presented more general polynomial-time attacks (with weaker leakage bounds) for the case where $e$ is not larger than $N^{1/2}$. A few years later, Blömer and May \cite{BM-partialkeyexposure} used heuristic Coppersmith techniques to obtain partial key exposure attacks for $e$ up to $N^{7/8}$. Let us recall their result and the proof idea.

\begin{lemma}[Theorem 11 in \cite{BM-partialkeyexposure}]
		For every $\varepsilon>0$ there exists $N_0$ such that for every $N\geq N_0$ the following holds:
		Let $(N,e)$ be a public RSA key with $\alpha=\log_N(e)\leq 7/8$, and define 
		\[
		\delta := \frac{1}{6} +\frac{1}{3}\sqrt{1+6\alpha}+\varepsilon.
		\]
		Assume that $d_0:=d\pmod M$ is known for some given $M \geq N^{\delta}$. Under Assumption \ref{resultant} explained in the proof sketch below, $N$ can be factored in polynomial time.
\label{thm:bmsmalle}
\end{lemma}

\begin{proof}[Sketch of the proof.]
	From $ed-1 = k\varphi(N) = k(N-(p+q-1))$, we obtain the bivariate congruence
	\[
	f(y,z)=y(N-z)-ed_0+1\equiv 0 \pmod{eM},
	\]
	with small root $(y_0,z_0) = (k,p+q-1)$. One easily observes that we may bound the root by
	\[
	|y_0| \leq Y=N^\alpha,
	\qquad 
	|z_0| \leq Z=3N^{1/2}.
	\]	
	Fix parameters $m,t$. We construct polynomial shifts which all vanish at $(y_0,z_0)$ modulo $(eM)^m$:
	\[
	g_{i,j}(y,z) = y^j(eM)^i f(y,z)^{m-i},
	\qquad
	i=0,\ldots,m,\quad j=0,\ldots,i,
	\]
	and
	\[
	h_{i,j}(y,z) = z^j(eM)^i f(y,z)^{m-i},
	\qquad
	i=0,\ldots,m,\quad j=1,\ldots,t.
	\]
	After the scaling $y \mapsto yY, z \mapsto zZ$, the coefficient vectors of these shifted polynomials form a lattice. Applying LLL (\Cref{thm:lll}) yields two short vectors in the lattice. For the specified choice of $\delta$ and suitable parameters $m$ and $t$, one can prove (\cite[Lemma 12]{BM-partialkeyexposure}) that these vectors are short enough to satisfy the conditions of Howgrave-Graham (\Cref{thm:hg}). 
	
	We obtain two corresponding polynomials $f_1$, $f_2$ that both vanish at $(y_0,z_0)$ over $\mathbb{Z}$, not only modulo $(eM)^m$. We find their common root by computing their resultant $\operatorname{Res}_y(f_1,f_2)$. However, this will only work if the resultant is nonzero.
	\begin{assumption} \label{resultant}
	 	The polynomials \(f_1,f_2\in\mathbb{Z}[y,z]\) obtained from this construction satisfy $ \operatorname{Res}_y(f_1,f_2)\not\equiv 0. $
	\end{assumption}
	Assuming Assumption \ref{resultant}, one may conclude that $(y_0,z_0)$ can be recovered and that $N$ can be factored.
\end{proof}

It is important to note that the modulus $M$ in Blömer-May's attack does not necessarily need to be a power of two. In our following result, we will make use of this fact by building $M$ as a product of a power of two (for known LSBs of $d$) and a primorial-like factor $P$ similar to the one used in \Cref{thm:sievebound}. This will allow us to run Blömer-May's attack even in cases where the known proportion of LSBs of $d$ is significantly below the threshold $\delta$. Our attack is similar to enumerations of unknown bits of $d$, but with an asymptotically significant reduction of the search space.

\begin{theorem}
		For every $\varepsilon>0$ there exists $N_0$ such that for every $N\geq N_0$ the following holds:
		Let $(N,e)$ be a public RSA key with prime $e$ and $\log_N(e)\leq 7/8$, and define $\delta$ as in \Cref{thm:bmsmalle}. Assume that $\ell$ LSBs $d\pmod{2^\ell}$ of $d$ are known. Let $\Delta:=N^{\delta}/2^\ell$ and $C:=(1+o(1))\log 2$. Under Assumption \ref{resultant}, $N$ can be factored in
		\[
		O\left(\Delta\cdot \exp\left(-C\log\Delta/\log\log\Delta\right)\cdot \polylog N\right)
		\]
		operations.
\label{thm:smalle}
\end{theorem}

\begin{proof} 
	Compute $P=\prod_{2<r\leq B} r$ for primes $r\neq e$ with $B$ minimal such that $P>\Delta$. 
	Let $M=2^\ell P>N^\delta$ and $d_{2^\ell}=d\pmod{2^\ell}$ be known, and set $v:=2^{-\ell}\pmod P$. 
	According to the Chinese Remainder Theorem, we then have
	\[
	d_M:=d\pmod M = d_{2^\ell}+2^\ell(v(d_P-d_{2^\ell}) \pmod P),
	\]
	where $d_P:=d\pmod P = e^{-1}(1+k\phi(N))\pmod P$.
	
	Let $\eta_P:=e^{-1}\pmod P$ and $\phi_P:=\phi(N)\pmod P$, write $d=\tilde{d}M+d_M$ and note that 
	\[
	e(\tilde{d}M+d_M) = e(\tilde{d}M +d_{2^\ell}+2^\ell(v(d_P-d_{2^\ell}) \pmod P))= 1+k(N-p-q+1),
	\] 
	hence
	\[
		 e2^\ell(v(d_P-d_{2^\ell}) \pmod P) = 1+k(N-p-q+1)-ed_{2^\ell}-eM\tilde{d}.
	\]	
	Considering the equation modulo $eM$, we get
	\begin{align*}
				 e2^\ell v(\eta_P(1+k\phi_P))&\equiv  1+k(N-p-q+1)+e2^\ell vd_{2^\ell}-ed_{2^\ell} \mod eM\\
				 ke2^\ell v\eta_P\phi_P &\equiv 1+k(N-p-q+1)+e2^\ell v(d_{2^\ell}-\eta_P )-ed_{2^\ell}\mod eM
	\end{align*}
	
	 We define the polynomial
	\[
	f_{\phi_P}(y,z):=y(N-e2^\ell v\eta_P\phi_P -  z) + e(2^\ell v(d_{2^\ell}-\eta_P )-d_{2^\ell}) +1
	\]	
	and note that $f_{\phi_P}(k,p+q-1)\equiv 0 \mod eM$. As a result, $f_{\phi_P}$ has a small root $(y_0,z_0)$ modulo $eM$ satisfying $y_0\leq N^{\alpha}$ and $z_0\leq 3N^{1/2}$. Assume we can compute all coefficients of $f_{\phi_P}$. We may then follow the approach by Blömer and May in \Cref{thm:bmsmalle} to retrieve $y_0$ and $z_0$ in polynomial time, using the same assumption concerning nonzero resultants.
	
	One observes that we can compute all coefficients of $f_{\phi_P}$ if we know $\phi_P$. We thus have to apply the procedure to all possible candidates for $\phi_P=\phi(N) \pmod P$, i.e., the elements of $\Phi_P$. The result then follows from \Cref{thm:sievebound} and \Cref{lem:phi}.
\end{proof}

\begin{remark}
	Since $\eta_P$ is the inverse of $e$ modulo $P$, the polynomial $f$ is only well-defined if $P$ is coprime to $e$. For simplicity, we assume $e$ to be prime and exclude it from the product $P$. While $e$ being prime is in line with standard practices and recommendations for the choice of $e$, it is stronger than what we actually need. We could also prove a version of Theorem \ref{thm:smalle} for composite $e$, but we would need to assume that its prime factorization does not contain an asymptotically significant number of those primes that we need for building the product $P$.
\end{remark}

\section{Attack on RSA with small private key}\label{sec:smalld}
Small private keys have long been known to make RSA vulnerable. The first major result in this direction was due to Wiener \cite{W-smalld} in 1990, who showed that $N$ can be factored efficiently if $d < N^{1/4}$ by recovering $d$ from the continued fraction expansion of $e/N$. Boneh and Durfee \cite{BD-smalld} later improved this bound using lattice methods, giving a heuristic attack for $d < N^{0.292}$. Ernst, Jochemsz, May and de Weger \cite{EJMW-partialkeyexposure} further developed these ideas in the setting of partial key exposure. In particular, they studied how much of a small private key $d=N^\beta$ needs to be known to enable factorization of $N$ in polynomial time. We give a reformulation of their result below.

\begin{lemma}[Theorem 3 in \cite{EJMW-partialkeyexposure}]
	For every $\varepsilon>0$ there exists $N_0$ such that for every $N\geq N_0$ the following holds:
	Let $(N,e)$ be a public RSA key, let $d=N^\beta$ be the private key, and define
	\[
	\delta := \frac{5}{6} -\frac{1}{3}\sqrt{1+6\beta}-\varepsilon.
	\]
	Assume that $d_0:=d\pmod M$ is known for some given $M \geq N^{\beta-\delta}$. Under Assumption \ref{resultant}, $N$ can be factored in polynomial time.
\label{thm:ernstsmalld}
\end{lemma}

\begin{proof}[Sketch of the proof.] 
	We consider $ed-1 = k\varphi(N) = k(N-(p+q-1))$ to obtain the bivariate congruence
	\[
	f(y,z)=y(N-z)-ed_0+1\equiv 0 \pmod{eM},
	\]
	with the same small root $(y_0,z_0) = (k,p+q-1)$ as in \Cref{thm:bmsmalle}. We may now bound this root by
	\[
	|y_0| \leq Y=N^\beta,
	\qquad 
	|z_0| \leq Z=3N^{1/2}.
	\]	
	The remainder of the proof follows the same approach as in \Cref{thm:bmsmalle}, with only a slight difference in notation. While $\delta$ in \Cref{thm:bmsmalle} denoted the size of the known part of $d$, here it denotes the size of the unknown part
	$
	N^\beta / M\leq N^{\delta}.
	$
\end{proof}

\begin{theorem}
	For every $\varepsilon>0$ there exists $N_0$ such that for every $N\geq N_0$ the following holds:
	Let $(N,e)$ be a public RSA key with prime $e$, let $d=N^{\beta}$ be the private key, and define $\delta$ as in \Cref{thm:ernstsmalld}. Assume that $\ell$ LSBs $d\pmod{2^\ell}$ of $d$ are known. Let $\Delta:=N^{\beta - \delta}/2^\ell$ and $C:=(1+o(1))\log 2$. Under Assumption \ref{resultant}, $N$ can be factored in
	\[
	O\left(\Delta\cdot \exp\left(-C\log\Delta/\log\log\Delta\right)\cdot \polylog N\right)
	\]
	operations.
\label{thm:smalld}
\end{theorem}

\begin{proof}
	Compute $P=\prod_{2<r\leq B} r$ for primes $r\neq e$ with $B$ minimal such that $P>\Delta$. 
	Let $M=2^\ell P>N^{\beta-\delta}$ and $d_{2^\ell}=d\pmod{2^\ell}$ be known. In addition, we define $v:=2^{-\ell}\pmod P$, $\eta_P:=e^{-1}\pmod P$ and $\phi_P:=\phi(N)\pmod P$.
	Consider the same polynomial as in the proof of \Cref{thm:smalle}, namely
	\[
	f_{\phi_P}(y,z):=y(N-e2^\ell v\eta_P\phi_P -  z) + e(2^\ell v(d_{2^\ell}-\eta_P )-d_{2^\ell}) +1,
	\]	
	and note that $f_{\phi_P}(k,p+q-1)\equiv 0 \mod eM$. As a result, $f_{\phi_P}$ has a small root $(y_0,z_0)$ modulo $eM$ satisfying $y_0\leq d < N^{\beta}$ and $z_0\leq 3N^{1/2}$. Under these conditions, together with Assumption \ref{resultant}, we have seen above that we may retrieve $y_0$ and $z_0$ in polynomial time. As in \Cref{thm:smalle}, the result follows from applying the procedure to all possible elements of $\Phi_P$.
\end{proof}

\section{Attack on RSA with small prime difference}\label{sec:smalldiff}
We recall that RSA is insecure if the prime factors $p$ and $q$ of $N$ are close together, meaning that $q-p$ is small. The reason for this is Fermat's factorization method, which works as follows. Consider
\[
(p-q)^2=(p+q)^2-4N=(p+q-2\sqrt{N})(p+q+2\sqrt{N}),
\]
which implies
\begin{equation}\label{eq:zbound}
0\leq p+q -2\sqrt{N}=\frac{(p-q)^2}{p+q+2\sqrt{N}}<\frac{(p-q)^2}{4\sqrt{N}}.
\end{equation}
For $L:=2\lceil \sqrt{N}\rceil$, we want to find a positive integer $z$ such that $L+z=p+q$. We can easily tell whether we have found the correct value $z$ by checking whether $(L+z)^2-4N$ is a square. The inequality above implies that, if $q-p \in O(N^{1/4})$, Fermat's method will factorize $N$ in $O(1)$ iterations. In general, the runtime cost consists of $O(\Lambda)$ square checks, where
\[
\Lambda:=\frac{(p-q)^2}{4\sqrt{N}}.
\]
An immediate improvement of this runtime complexity via the hyperbolic sieve is proved in \cite[Theorem 4.6]{Hit-subsetsum}, achieving a reduction of the number of required square checks to 
\[
\widetilde{O}\left(\Lambda\cdot \exp\left(-C\log\Lambda/\log\log\Lambda\right)\right)
\]
with $C=(1+o(1))\log 2$, while keeping the space complexity negligible.

\vspace{12pt}
We will now consider a substantial asymptotic improvement over Fermat's original method for factoring semiprimes with close factors. The attack has been mentioned in prior work (\cite{Hit-BSGS}, \cite{Hit-subsetsum}), but has never been stated and proved as a formal theorem. In addition, an important ingredient of the following rigorous proof only became available recently.
\begin{lemma}
	\label{lem:large_order}
	Let $N \geq 3$ and $D<N-1$ be positive integers.
	There is a deterministic algorithm that either finds a factor of $N$, or some $\alpha \in\Z_N^*$ such that $\ord_N(\alpha) > D$, and runs in time $\widetilde{O}(D^{1/2}\cdot \polylog(N))$.
\end{lemma}

\begin{proof}
	This is Theorem 1.1 in \cite{HH-largeorder}. For our application, it will be important that $D$ is allowed to be any positive integer up to $N-1$. Similar previous results assumed a lower bound for $D$ that was exponential in $\log N$. 
\end{proof}

To the best of our knowledge, the following result yields the fastest known way for factoring RSA moduli with small prime difference, both theoretically and in practice. 

\begin{theorem}
	\label{thm:smallpq}
Let $N=pq\in\N$ be semiprime, let $\Delta<N-1$ be a positive integer, and set $C:=(1+o(1))\log 2$. If
$
(q-p)^2/4\sqrt{N}\leq \Delta
$
then one can deterministically factor $N$ in
\[
O\left(\Delta^{1/2}\cdot \exp\left(-\tfrac{C}{2}\log \Delta/\log\log \Delta\right)\cdot \polylog N\right)
\]
bit operations. 
\end{theorem}

\begin{proof}
Let $S:=p+q$ and $z:=S-2\lceil\sqrt{N}\rceil$. From (\ref{eq:zbound}) it follows that $z\leq \Delta$. As a consequence of $\phi(N)=N-S+1$ we may write
\begin{equation}\label{eq:bsgs}
	\alpha^{z}\equiv \alpha^{N+1-2\lceil\sqrt{N}\rceil} \mod N
\end{equation}
for every $\alpha\in\Z$ coprime to $N$. We will use a babystep-giantstep method to solve the discrete logarithm problem (\ref{eq:bsgs}) for $z$, which allows us to compute $S$. Once we have found $S$, we can solve the quadratic equation $X^2-SX+N=0$ to obtain the factors $p$ and $q$ of $N$ as solutions.

We start by considering a different representation of $z$ and constructing the baby- and the giantsteps. Let $B\in\N$ and consider the product $P$ of all odd primes up to $B$, 
\[
P=\prod_{2<r\leq B} r.
\] 
Consider any decomposition $P=ML$ of $P$ and note that $\gcd(M,L)=1$. For $z_M := z \pmod M$ and $z_L := z\pmod L$, we may then write
$
z=(kL+j)M+z_M
$
with $j,k\in\N$ such that $j = M^{-1}(z_L-z_M)\pmod L$ and $k\leq \Delta/P$. One observes that we also have
\[
z=(\bar{k}L+\bar{j})M+z_M
\]
with $\bar{j}=M^{-1}z_L \pmod L - M^{-1}z_M\pmod L$. Note that either $\bar{j}=j$ or $\bar{j}=j-L$. In any case, we obtain $\bar{k}\leq \Delta/P+1$.
We may then rewrite (\ref{eq:bsgs}) to
\[
\alpha^{z_M + (M^{-1}z_M)\pmod L}\equiv \alpha^{-\bar{k}P}\cdot \alpha^{-M^{-1}z_L\pmod L}\cdot\alpha^{N+1-2\lceil\sqrt{N}\rceil}\mod N.
\]
The only unknown on the left-hand side (the ``babysteps'') is $z_M$, and the only unknowns on the right-hand side (the ``giantsteps'') are $z_L$ and $\bar{k}$. So the number of babysteps is equal to all possible  candidates for $z_M$, and the number of giantsteps is $z_L\cdot (\Delta/P+1)$.

We now choose $B=(1+o(1))\log\Delta$ such that $P\geq\Delta$ is minimal, implying that $P\in \widetilde{O}(\Delta)$. This implies that $\bar{k}\in\widetilde{O}(\log\Delta)$. In addition, \Cref{thm:sievebound} implies that the set $\mathcal{S}_P$ containing all residue candidates for $z\pmod P$ is of size
\[
\widetilde{O}\left(P\cdot \exp\left(-C\log P/\log \log P\right)\right) = \widetilde{O}\left(\Delta\cdot \exp\left(-C\log \Delta/\log\log \Delta\right)\right).
\]
One easily observes that $\mathcal{S}_P$ is the image of $\mathcal{S}_M \times \mathcal{S}_L$ under the inverse Chinese Remainder Theorem isomorphism. As a consequence, $|\mathcal{S}_P|=|\mathcal{S}_M|\cdot |\mathcal{S}_L|$. We want to find a decomposition of $P=ML$ such that $\mathcal{S}_M$ and $\mathcal{S}_L$ are about the same size, namely both around
\[
\widetilde{O}\left(\Delta^{1/2}\cdot \exp\left(-\tfrac{C}{2}\log \Delta/\log \log \Delta\right)\right).
\]
Considering our estimates above, we should choose $M$ as the product of odd primes such that $M\in\widetilde{O}(P^{1/2})$, and $L=P/M$. 

The algorithm now works as follows: We start by applying Lemma \ref{lem:large_order}. The cost is $\widetilde{O}(L^{1/2}\cdot \polylog N)\subseteq \widetilde{O}(\Delta^{1/4}\cdot \polylog N)$, which is negligible. We either factor $N$ or we find an element $\alpha\in\Z_N^*$ such that $\ord_N(\alpha)>L$. Assume the latter. For all $b\in\mathcal{S}_M$, we then compute the babysteps
\[
\mathcal{B}(b):=\alpha^{b + (M^{-1}b)\pmod L}\pmod N.
\]
For all $g\in\mathcal{S}_L$ and $k<\Delta/P+1\in \widetilde{O}(\log\Delta)$, we compute the giantsteps
\[
\mathcal{G}(g,k):=\alpha^{-\bar{k}P}\cdot \alpha^{-M^{-1}g\pmod L}\cdot\alpha^{N+1-2\lceil\sqrt{N}\rceil}\pmod N.
\]
Using standard techniques for modular exponentiation, this can be done within the claimed runtime. Our discussion above now implies that there exist $z_M$, $z_L$ and $\bar{k}$ such that 
\[
\mathcal{B}(z_M)\equiv \mathcal{G}(z_L,\bar{k}) \mod N.
\]
Such matches can be found by applying merge sort (\cite[Chap. 2.2]{SW-algorithms}) to the lists of the baby- and the giantsteps and performing a subsequent scan for matches. Since $\ord_N(\alpha)>L$, the $|\mathcal{S}_L|$ elements $\alpha^{-M^{-1}g\pmod L} \pmod N$ on the giantstep-side are all pairwise distinct. Together with the fact that $\bar{k}\in\widetilde{O}(\log \Delta)$, this implies that the total number of matches between the baby- and the giantsteps is bounded by $ \widetilde{O}(|\mathcal{S}_M|)$. The overall cost for finding all these matches is thus within the claimed runtime (\cite[Prop. F on p. 272]{SW-algorithms}).

For each found match $\mathcal{B}(b)\equiv \mathcal{G}(g,k)\mod N$, we compute the resulting candidate for $z=z\pmod P$ from the values $b\pmod M$ and $g\pmod L$. For the correct match $b=z_M$ and $g=z_L$, using $S=z+2\lceil\sqrt{N}\rceil$ as coefficient in the quadratic equation mentioned in the beginning of the proof will allow us to find the factor $N$. Since each such equation can be solved in $O(\polylog N)$ time, this concludes the proof.
\end{proof}

\begin{corollary}
	Let $N=pq\in\N$ be semiprime, let
	$
	\Lambda:=(q-p)^2/4\sqrt{N}
	$
	and $C:=(1+o(1))\log 2$. Then one can deterministically factor $N$ in
	\[
	O\left(\Lambda^{1/2}\cdot \exp\left(-\tfrac{C}{2}\log \Lambda/\log\log \Lambda\right)\cdot \polylog N\right)
	\]
	bit operations. 
\end{corollary}

\begin{proof}
	Apply Theorem \ref{thm:smallpq} with $\Delta=2^i$ for $i=1,2,3,\ldots$ until the factorization of $N$ is found. We succeed as soon as $\Delta>\Lambda$, which will be the case after at most $\log_2 \Lambda \in O(\polylog N)$ iterations.
\end{proof}

\bibliographystyle{amsplain}
\bibliography{references}

\end{document}